\documentclass[8pt]{article}
\usepackage{graphicx}
\usepackage{pgfplots}
\usepackage{newlfont}
\usepackage[colorlinks,linkcolor=blue,citecolor=blue,anchorcolor=blue,bookmarksopen,pdfpagetransition={Wipe}]{hyperref}
\usepackage{amsmath,amsthm,amssymb}
\usepackage{expl3}
\usepackage{amsfonts}
\newtheorem{thm}{Theorem}[section]

\newtheorem{lemma}[thm]{Lemma}

\newtheorem{definition}[thm]{Definition}

\newtheorem{example}[thm]{Example}
\numberwithin{equation}{section}

\usepackage{color} 

\begin{document}
\title{\textbf{High-order numerical method for two-dimensional Riesz space fractional advection-dispersion equation}}
\author{A. Borhanifar$^a$\footnote{borhani@uma.ac.ir}\
, M. A. Ragusa$^b$$^c$\footnote{Corresponding author: maragusa@dmi.unict.it} and S. Valizadeh$^a$\footnote{valizadeh.s@uma.ac.ir}}

\date{}
\maketitle
\begin{center}
$^a$Department of Mathematics, University of Mohaghegh Ardabili,
56199-11367 Ardabil, Iran\\
$^b$Dipartimento di Matematica e Informatica, Universit\`{a} di Catania, Viale Andrea Doria,
 6-95125 Catania, Italy\\
$^c$RUDN University, 6 Miklukho - Maklay St, Moscow, 117198, Russia
\end{center}
\begin{abstract}
\indent In this paper, by combining of fractional centered difference approach with alternating direction implicit method, we introduce a mixed difference method for solving two-dimensional Riesz space fractional advection-dispersion equation. The proposed method is a fourth order centered difference operator in spatial directions and second order Crank-Nicolson method in temporal direction. By reviewing the consistency and stability of the method, the convergence of the proposed method is achieved. Several numerical examples are considered aiming to demonstrate the validity and applicability of the proposed technique.
\end{abstract}

\vskip .3cm

\section{Introduction}
\indent \hskip .65cm

Many researchers focused on fractional partial differential equations due to their useful applications in many real-world models, modeling with the least error and overlapping physically with scientific issues.

Fractional partial differential equations are mainly classified to the time, space, and time-space fractional partial differential equations. Of all these types, space fractional partial differential equations are containing fractional diffusion equation, fractional Fokker-Planck equation, fractional advection-dispersion equation (see e.g. \cite{Adams19923293,Baeumer20011543,Gala20121741,Gala20161271}).
The fractional advection-dispersion equation provides a useful description of chemical and contaminant transport in heterogeneous aquifers \cite{Adams19923293,Baeumer20011543},
abnormal mass absorption in solids \cite{Povstenko2017203} and densities of plumes in spread proportionally to time-dependent of fractional order \cite{Benson20001403}.
Although there is a plenty of research in one dimensional Riesz space fractional advection-dispersion equation (for example, see \cite{Abdi-mazraeh2015818,Ding2017759,Ding2015218,Irandout-pakchin2014913,Lakestani20121149,Manafian2017243,Micu20061950,Popolizio20131975,
Rahman2014264,Shen2008850,Valizadeh2020,Zhang2014266}), there are few works in the two dimensional case \cite{Chen2008295,Zeng20142599,Zhang20172063}. In this research, we numerically solve the general two dimensional Riesz space fractional advection-dispersion equation using modified Crank-Nicolson Alternating Direction Implicit (ADI) method. The outline of the paper is as follows. Problem description and its applications are presented in the next section. Section 3 is devoted to basic definitions and preliminaries. Section 4 is dedicated to the construction and explanation of the numerical scheme. In Section 5, consistency and stability of the numerical scheme are theoretically analyzed. Numerical examples are given in Section 6. In Section 7 the conclusions of the obtained results are given.

\section{Problem description}

Let $\Pi$ be a rectangular domain in $\mathbb{R}^2$ with boundary $\Upsilon=\partial \Pi$ and $\mathfrak{T}=(0,T]$ be the time interval,
$T > 0$. Consider the following two dimensional Riesz space fractional partial differential equation for a solute concentration of material $u$

\begin{eqnarray*}
\frac{\partial u(x,y,t)}{\partial
t}=d_{\alpha}\frac{\partial^{\alpha}u(x,y,t)}{\partial
|x|^{\alpha}}+c_{\beta}\frac{\partial^{\beta}u(x,y,t)}{\partial
|x|^{\beta}}+d_{\mu}\frac{\partial^{\mu}u(x,y,t)}{\partial
|y|^{\mu}}+c_{\nu}\frac{\partial^{\nu}u(x,y,t)}{\partial
|y|^{\nu}}
\end{eqnarray*}
\begin{eqnarray}\label{math201803252034}
+s(x,y,t),\quad (x,y,t)\in\Pi\times \mathfrak{T},
\end{eqnarray}
\begin{eqnarray}\label{math201803252035}
u(x,y,0)=f(x,y), \quad (x,y)\in\Pi,
\end{eqnarray}
\begin{eqnarray}\label{math201803252036}
u(x,y,t)=0, \quad (x,y,t)\in \Upsilon\times \mathfrak{T},
\end{eqnarray}

where $f:\Pi \rightarrow \mathbb{R}$ is a known function, $c_{\beta}$ and $c_{\nu}$ are the average fluid velocities in the x- and y-directions, while $d_{\alpha}$ and $d_{\mu}$ are the dispersion coefficients in the x- and y-directions, respectively. We restrict ourselves to $1 < \alpha, \mu \leq 2$, $0 < \beta, \nu < 1$ and assume that $d_{\alpha},d_{\mu}>0$ and $c_{\beta},c_{\nu}\geq 0$. The solution $u=u(x,y,t)$ is assumed to be sufficiently smooth and has the necessary continuous partial derivatives up to certain orders. The specific cases for equation (\ref{math201803252034}) is discussed in the years 2007, 2012 and 2014 by Valdes-Parada et al. \cite{Valdes-Parada2007}, Ding and Zhang \cite{Ding20121135} and Zeng et al. \cite{Zeng20142599} as follows:
\begin{itemize}
\item[(i)] \emph{The chemical reaction occurs only at the surface of the solid phase} \cite{Valdes-Parada2007}. The
conservation equation that governs the transport process in the $\gamma$-phase is introduced by
\begin{eqnarray}\label{rev3-2}
\frac{\partial c_{\gamma}}{\partial
t}+\nabla \bullet N_{\gamma}=0\quad \mbox{in the}\quad \gamma\mbox{-phase},
\end{eqnarray}
where $c_{\gamma}$ and $N_{\gamma}$ are the concentration and the molar flux of the reactant in the $\gamma$-phase, respectively.\\
A vector version of the fractional Fick's equation is given by:
\begin{eqnarray}\label{rev3-3}
N_{\gamma}=-\emph{D}_{\alpha,\gamma}\nabla^{\alpha}c_{\gamma},\quad \alpha \in (0,1],
\end{eqnarray}
where $\emph{D}_{\alpha,\gamma}$ is the fractional diffusion coefficient, and the fractional derivative operator is defined in the Riemann-Liouville's sense.
The combination of the two Eqs. (\ref{rev3-2}) and (\ref{rev3-3}) yields
\begin{eqnarray}\label{rev3-4}
\frac{\partial c_{\gamma}}{\partial
t}=\nabla \bullet(\emph{D}_{\alpha,\gamma}\nabla^{\alpha}c_{\gamma}),\quad \mbox{in the}\quad \gamma\mbox{-phase}.
\end{eqnarray}

\item[(ii)] \emph{1D Riesz space fractional advection-dispersion equation} \cite{Ding20121135}. Ding and Zhang considered in
\begin{eqnarray}\label{rev3-6}
\frac{\partial u(x,t)}{\partial
t}=K_{\alpha}\frac{\partial^{\alpha}u(x,t)}{\partial
|x|^{\alpha}}+K_{\beta}\frac{\partial^{\beta}u(x,t)}{\partial
|x|^{\beta}}
\end{eqnarray}
where the Riesz fractional derivatives are used, $u$ is a solute concentration; $K_{\alpha}$ and $K_{\beta}$ represent the dispersion coefficient and the average fluid velocity. They restrict $1 < \alpha\leq 2$, $0 < \beta< 1$, and assume that $K_{\alpha}>0$ and $K_{\beta}\geq 0$.
This is a special case for (\ref{math201803252034}) where $d_{\alpha} \equiv K_{\alpha}$, $c_{\beta} \equiv K_{\beta}$, $d_{\mu}=c_{\nu}=0$ and $s(x,y,t)=0$.

\item[(iii)] \emph{2D nonlinear chemical reaction} \cite{Zeng20142599}. By briefly typing sentences for fractional derivatives in (\ref{rev3-4}) and apply it to another dimension and also getting non-linear term,  two-dimensional Riesz space fractional nonlinear reaction-diffusion equation is introduced as
\begin{eqnarray}\label{rev3-5}
\partial_{t} u=K_{x}\frac{\partial^{2\alpha_{1}}u}{\partial
|x|^{2\alpha_{1}}}+K_{y}\frac{\partial^{2\alpha_{2}}u}{\partial
|y|^{2\alpha_{2}}}+F(u)+f(x,y,t)
\end{eqnarray}
in which $\frac{1}{2}<\alpha_{1},\alpha_{1}<1$, $K_{x},K_{y}>0$. This is special form of (\ref{math201803252034}) for $c_{\beta}=c_{\nu}=F(u)=0$.
\end{itemize}
The equation that is considered in this paper relates to a modeling of the fractional Fick's law in porous media in both x- and y- directions. For this reason the basic definitions and lemmas will be considered in the next section.

\section{Preliminaries}

In this section, we consider some important definitions, lemmas and notations which are useful for the further progresses in this paper.

\begin{definition}
The left-and right-sided Riemann-Liouville fractional derivatives of order $\gamma$
of $f(x)$ that be a continuous and Smooth function of order $n$ on $[0,L]$ are defined respectively as \cite{Podlubny1999},
\begin{eqnarray}\label{math21}
_{0} D_{x}^{\gamma}f(x)=\frac{1}{\Gamma(n-\gamma)}\frac{d^{n}}{d x^{n}}\int_{0}^{x}\frac{f(\xi)}{(x-\xi)^{\gamma-n+1}}d\xi,
\end{eqnarray}
\begin{eqnarray}\label{math22}
_{x} D_{L}^{\gamma}f(x)=\frac{(-1)^n}{\Gamma(n-\gamma)}\frac{d^{n}}{d x^{n}}\int_{x}^{L}\frac{f(\xi)}{(\xi-x)^{\gamma-n+1}}d\xi,
\end{eqnarray}
where $\Gamma(\cdot)$ is the Euler gamma function, $n-1 < \gamma \leq n$, $n\in N$ and $n$ is the smallest integer greater than $\gamma$.
\end{definition}

\begin{definition}
The Riesz fractional derivatives of order $\gamma$ of $f(x)$ function on $[0,L]$ is defined as \cite{Gorenflo1999231},
\begin{eqnarray}\label{math21}
\frac{\partial^{\gamma}f(x)}{\partial
|x|^{\gamma}}=-\kappa_{\gamma}\{_{0}D_{x}^{\gamma}f(x)+_{x}D_{L}^{\gamma}f(x)\},
\end{eqnarray}
where
$\kappa_{\gamma}=\frac{1}{2cos(\frac{\pi\gamma}{2})}$, $n-1< \gamma \leq n$ and $\gamma\neq1$.
\end{definition}

\begin{lemma}\label{math201803072224}
\cite{Meyer2000} The eigenvalues of the following tridiagonal Toeplitz matrix
\begin{equation*}\label{3.42}
A=\left(
\begin{array}{cccccccc}
  b & a &   &   &   &   &  \\
  c & b & a &   &   & O &   \\
    & c & \cdot & \cdot &   &   & \\
    &   & \cdot & \cdot & \cdot &   & \\
    &   &   & \cdot & \cdot & \cdot & \\
    & O &   &   & \cdot & \cdot & a   \\
    &   &   &   &   & c & b
\end{array}
\right)_{n\times n}
\end{equation*}
are given by
\begin{equation}\label{3.41}
\lambda_{j}=b+2a\sqrt{c/a}\cos(j\pi/(n+1)), \quad j=1,2,...,n,
\end{equation}
while the corresponding eigenvectors are:
\begin{equation*}
x_{j}=\left(
\begin{array}{c}
  (c/a)^{1/2}\sin(1j\pi/(n+1)) \\
  (c/a)^{2/2}\sin(2j\pi/(n+1)) \\
  (c/a)^{3/2}\sin(3j\pi/(n+1)) \\
  \vdots \\
  (c/a)^{n/2}\sin(nj\pi/(n+1)) \\
\end{array}
\right), \quad j=1,2,...,n,
\end{equation*}
i.e., $A x_{j}=\lambda_{j} x_{j}$, $j=1,2,...,n$. Moreover, the matrix $A$ is diagonalizable and $P=(x_{1}\quad x_{2}\quad ...\quad x_{n})$ diagonalizes $A$,
i.e., $P^{-1}AP=D$, where $D=diag(\lambda_{1}\quad \lambda_{2}\quad ...\quad \lambda_{n})$.
\end{lemma}

Here we consider the approximation with step $\mathit{h}$ of the Riesz fractional derivative that obtained by calculating the appropriate coefficients for the fractional central difference by applying Fourier transform \cite{Ding2015218}\\
\begin{eqnarray}\label{math4}
\frac{\partial^{\gamma}u(x)}{\partial
|x|^{\gamma}}\approx-h^{-\gamma}\sum_{r=-1}^{1}\varrho_{r}^{(\gamma)}\mathcal{H}_{r}^{(\gamma)}u(x),\quad 0< \gamma \leq 2 \quad and \quad \gamma\neq1,
\end{eqnarray}
where
\begin{eqnarray*}\label{math4-4}
\varrho_{-1}^{(\gamma)}=\varrho_{1}^{(\gamma)}=-\frac{\gamma}{24}, \quad \varrho_{0}^{(\gamma)}=\frac{\gamma}{12}+1,
\end{eqnarray*}
\begin{eqnarray}\label{math4-1}
\mathcal{H}_{r}^{(\gamma)}u(x)=\sum_{k=-\infty}^{\infty}\omega_{k}^{(\gamma)}u(x-(k+r)h),
\end{eqnarray}
and all coefficients $\omega_{k}^{(\gamma)}$ are defined by
\begin{eqnarray}\label{math4-2}
\omega_{k}^{(\gamma)}=\frac{(-1)^{k}\Gamma(\gamma+1)}{\Gamma(\frac{\gamma}{2}-k+1)\Gamma(\frac{\gamma}{2}+k+1)}, \quad k=0,\pm 1,\pm 2,...,
\end{eqnarray}

By inserting the values of the $\varrho_{r}^{(\gamma)}$ and series in formula (\ref{math4-1}) into the formula (\ref{math4}), the corresponding operator of the Riesz fractional derivative of order $\gamma$ that we denote via $\mathfrak{D}_{\gamma,x}$ will be as follows
\begin{eqnarray*}
\mathfrak{D}_{\gamma,x}u(x)=-\frac{\gamma}{24}\sum_{k=-\infty}^{\infty}\frac{\omega_{k}^{(\gamma)}u(x-(k-1)h)}{h^{\gamma}}
+(\frac{\gamma}{12}+1)\sum_{k=-\infty}^{\infty}\frac{\omega_{k}^{(\gamma)}u(x-k h)}{h^{\gamma}}
\end{eqnarray*}
\begin{eqnarray*}
-\frac{\gamma}{24}\sum_{k=-\infty}^{\infty}\frac{\omega_{k}^{(\gamma)}u(x-(k+1)h)}{h^{\gamma}}
\end{eqnarray*}

We survey the properties of the coefficients $\omega_{k}^{(\gamma)}$ that are appearing at the approximate formula for Riesz fractional derivatives.

\begin{lemma}\label{Thm99}
    \cite{Celik20121743} The coefficients  $\omega_{k}^{(\gamma)}$ for $k\in \mathbb{Z}$ in (\ref{math4-2}) satisfy:\\
(a) $\omega_{0}^{(\gamma)}\geq0$, $\omega_{-k}^{(\gamma)}=\omega_{k}^{(\gamma)}\leq 0$ for all $\mid k \mid \geq 1$,\\
(b) $\sum_{k=-\infty}^{\infty}\omega_{k}^{(\gamma)}=0$,\\
(c) For any positive integer $n$ and $m$ with $n<m$, we have $\sum_{k=-m+n}^{n}\omega_{k}^{(\gamma)}>0$,\\
(d) $\mid 2\sin(\frac{z}{2}) \mid^{\gamma}=\sum_{k=-\infty}^{\infty}\omega_{k}^{(\gamma)}e^{-ik z}$.
\end{lemma}

\begin{thm}\label{math201802221502}
Let $f\in C^{7}(\mathbb{R})$ and all derivatives up to order seven belong to $L_{1}(\mathbb{R})$. Then
\begin{eqnarray}\label{math201802221503}
\frac{\partial^{\gamma}f(x)}{\partial
|x|^{\gamma}}=-h^{-\gamma}\sum_{r=-1}^{1}\varrho_{r}^{(\gamma)}\mathcal{H}_{r}^{(\gamma)}f(x)+\mathcal{O}(h^{4}),
\end{eqnarray}
when $h \rightarrow 0$ and $\frac{\partial^{\gamma}f(x)}{\partial
|x|^{\gamma}}$ is the Riesz fractional derivative for $\gamma \in(0,1)\cup(1,2]$.
\end{thm}

\begin{proof}
Let
\begin{eqnarray}\label{math201802271726}
\varepsilon(x,h)=\frac{\partial^{\gamma}f(x)}{\partial
|x|^{\gamma}}+h^{-\gamma}\sum_{r=-1}^{1}\varrho_{r}^{(\gamma)}\mathcal{H}_{r}^{(\gamma)}f(x).
\end{eqnarray}
Applying the Fourier transformation defined as $\widehat{f}(\xi)=\mathcal{F}\{f(x)\}=\int_{-\infty}^{\infty}e^{-i\xi x}f(x) dx$, $\xi \in \mathbb{R}$ to the Eq. (\ref{math201802271726}) yields

\begin{eqnarray*}
\widehat{\varepsilon}(\xi,h)=-\mid \xi \mid^{\gamma}\widehat{f}(\xi)-[\frac{\gamma}{24}\sum_{k=-\infty}^{\infty}\frac{\omega_{k}^{(\gamma)}e^{-i(k-1)\xi h}}{h^{\gamma}}
-(\frac{\gamma}{12}+1)\sum_{k=-\infty}^{\infty}\frac{\omega_{k}^{(\gamma)}e^{-ik\xi h}}{h^{\gamma}}
\end{eqnarray*}
\begin{eqnarray}\label{math201802280048}
+\frac{\gamma}{24}\sum_{k=-\infty}^{\infty}\frac{\omega_{k}^{(\gamma)}e^{-i(k+1)\xi h}}{h^{\gamma}}]\widehat{f}(\xi),
\end{eqnarray}

According to being $\mid \frac{2\sin(\frac{\xi h}{2})}{h} \mid$ the generating function of the coefficients $\omega_{k}^{(\gamma)}$ (formula in Lemma \ref{Thm99} part (d) ), the above relation will be as follows

\begin{eqnarray}\label{math201802222020}
\widehat{\varepsilon}(\xi,h)=[-\mid \xi \mid^{\gamma}-(\frac{\gamma}{12}\cos{\xi h}-(\frac{\gamma}{12}+1))\mid \frac{2\sin(\frac{\xi h}{2})}{h} \mid^{\gamma}]\widehat{f}(\xi),
\end{eqnarray}
also, we have the following relation based on Taylor's expansion of the function $\sin(\frac{\xi h}{2})$ centered at $\xi=0$
\begin{eqnarray*}
\mid \frac{2\sin(\frac{\xi h}{2})}{h} \mid^{\gamma}=\mid \xi \mid^{\gamma}[1-\frac{\gamma}{24}(\xi h)^{2}+(\frac{1}{1920}+\frac{\gamma-1}{1152})\gamma(\xi h)^{4}-(\frac{1}{322560}+\frac{\gamma-1}{46080}
\end{eqnarray*}
\begin{eqnarray}\label{math201802222022}
+\frac{(\gamma-1)(\gamma-2)}{82944})\gamma(\xi h)^{6}+ \mathcal{O}((\xi h)^{8})].
\end{eqnarray}
Putting formula (\ref{math201802222022}) in formula (\ref{math201802222020}) gives the following relation:
\begin{eqnarray}\label{math201802272214}
\widehat{\varepsilon}(\xi,h)=\mid \xi \mid^{\gamma}[\frac{\gamma}{288}-\frac{\gamma^{2}}{576}+(\frac{1}{1920}+\frac{\gamma-1}{1152})\gamma](\xi h)^{4}\widehat{f}(\xi)
+ \mathcal{O}((\xi h)^{6}).
\end{eqnarray}
Since $f\in C^{7}(\mathbb{R})$ and all derivatives up to order seven belong to $L_{1}(\mathbb{R})$, there exists a positive constant $C_{0}$ that
\begin{eqnarray}\label{math201802272314}
\mid \widehat{f}(\xi) \mid \leq C_{0}(1+\mid \xi \mid)^{-7}.
\end{eqnarray}
Therefore, we have from (\ref{math201802272214}) and (\ref{math201802272314})
\begin{eqnarray*}\label{math201802272318}
\mid \widehat{\varepsilon}(\xi,h) \mid \leq \mathbb{C}h^{4}\mid \xi \mid^{4+\gamma}C_{0}(1+\mid \xi \mid)^{-7}\leq C h^{4}(1+\mid \xi \mid)^{4+\gamma}(1+\mid \xi \mid)^{-7}
\end{eqnarray*}
\begin{eqnarray}\label{math201802272318}
=C h^{4}(1+\mid \xi \mid)^{\gamma-3}.
\end{eqnarray}
where $C=\mathbb{C}C_{0}$ is independent of $\xi$. That is, the inverse Fourier transform of the function $\widehat{\varepsilon}(\xi,h)$ exists for $\gamma \in(0,1)\cup(1,2]$.

Therefore, taking the inverse Fourier transform in both sides of (\ref{math201802280048}) and using (\ref{math201802272318}) gives

\begin{eqnarray*}
\mid \varepsilon(x,h)\mid
\leq \frac{1}{2 \pi} \int_{-\infty}^{\infty}\mid\widehat{\varepsilon}(\xi,h)\mid d\xi \leq \int_{-\infty}^{\infty}C h^{4}(1+\mid \xi \mid)^{\gamma-3}d\xi=C'h^{4},
\end{eqnarray*}
where $C'=\frac{C}{(2-\gamma)\pi}$.\\
Hence
\begin{eqnarray*}
\varepsilon(x,h)=\mathcal{O}(h^{4}).
\end{eqnarray*}
\end{proof}

\begin{lemma}\label{math201803072254}
The matrix form of the operator $\mathfrak{D}_{\gamma,x}$ is symmetric positive definite.
\end{lemma}

\begin{proof}
Without loss of generality, suppose to be $m$ the number of nodal points for the Riesz fractional derivative approximation. The matrix form of the operator $\mathfrak{D}_{\gamma,x}$ is as follow
\begin{eqnarray*}
M(\mathfrak{D}_{\gamma,x})=A_{x}^{\gamma}\times B_{x}^{\gamma},
\end{eqnarray*}
where $A_{x}^{\gamma}$ and $B_{x}^{\gamma}$ are two symmetric matrices of order $m$, which have the entries $\varrho_{r}^{(\gamma)}$ and $\omega_{k}^{(\gamma)}$, respectively,
\begin{eqnarray*}
A_{x}^{\gamma}=tridiag\{\varrho_{-1}^{(\gamma)},\varrho_{0}^{(\gamma)},\varrho_{1}^{(\gamma)}\}=
tridiag\{-\frac{\gamma}{24},1+\frac{\gamma}{12},-\frac{\gamma}{24}\},
\end{eqnarray*}
i.e., the matrix $A_{x}^{\gamma}$ is the submatrix of matrix $A$ of order $m$, considered in the Lemma \ref{math201803072224}, for $a=c=-\frac{\gamma}{24}$ and $b=1+\frac{\gamma}{12}$.
\begin{eqnarray*}
(B_{x}^{\gamma})_{i,j}=\omega_{\mid i-j\mid}^{(\gamma)}, \quad i,j=1,2,...,m
\end{eqnarray*}
Based on the Lemma \ref{math201803072224} eigenvalues and eigenvectors of matrix $A_{x}^{\gamma}$ are $\mu_{k}=1+\frac{\gamma}{12}\sin^{2}\frac{k\pi}{2(m+1)}$ and $y_{k}=(\sin(\frac{k \pi}{m+1})\quad \sin(\frac{2 k \pi}{m+1})\quad ...\quad\sin(\frac{m k \pi}{m+1}))^{T}$ for $k=1,2,...,m$ and  also $A_{x}^{\gamma}$ is a real diagonally dominant matrix and according to the Lemma \ref{Thm99} eigenvalues of matrix $B_{x}^{\gamma}$ are positive and $B_{x}^{\gamma}$ is a real diagonally dominant matrix and obviously the matrices $A_{x}^{\gamma}$ and $B_{x}^{\gamma}$ are symmetric. Therefore the matrices $A_{x}^{\gamma}$ and $B_{x}^{\gamma}$ are symmetric positive definite. The proof of the lemma is finished.
\end{proof}

\section{Implementation high order ADI scheme}

In this section, we apply modified fractional centered difference ADI method to solve two dimensional RSFADE.
For the numerical approximation scheme, let $\Pi$ be a finite domain satisfying $\Pi=[x_{L},x_{R}]\times[y_{L},y_{R}]$ and $0\leq t \leq T$, we introduce a uniform grid of mesh points $(x_{i},y_{j},t_{s})$, $\Delta x=\frac{x_{R}-x_{L}}{m_{1}}$ is the spatial grid size in $x$-direction, with partition $x_{i} = x_{L} + i\Delta x$ for $i=0,1,2,...,m_{1}$; $\Delta y=\frac{y_{R}-y_{L}}{m_{2}}$ is the spatial grid size in $y$-direction, with partition $y_{j} = y_{L} + j\Delta y$ for $j=0,1,2,...,m_{2}$ and $\Delta t=\frac{T}{N}$ is time step, with partition $t_{n}=n\Delta t$ for $n=0,1,...,N$, where $m_{1}$, $m_{2}$, and $N$ are being positive integers.
For any function $\phi(x, y, t)$, we let $\phi_{i,j}(t)=\phi(x_{i},y_{j},t)$ and $\phi_{i,j}^{n}=\phi(x_{i},y_{j},t_{n})$.\\


Firstly, we use modified fractional centered difference scheme to discrete Riesz derivatives with the subject of that these difference schemes operate on functions belong to $C^{7}(\mathbb{R})$.
The following discretization formulas are achieved based on the expression in Theorem \ref{math201802221502}.

\begin{eqnarray*} \label{math201804010425}
\frac{\partial^{\alpha}u(x_{i},y,t)}{\partial
|x|^{\alpha}}=-\Delta x^{-\alpha}\sum_{r=-1}^{1}\varrho_{r}^{(\alpha)}\sum_{k=-\infty}^{\infty}\omega_{k}^{(\alpha)}u(x_{i-k-r},y,t)
\end{eqnarray*}
\begin{eqnarray} \label{math201804010425}
+\mathcal{O}(\Delta x^{4}),\quad 1< \alpha \leq2 ,
\end{eqnarray}

\begin{eqnarray*} \label{math201804010426}
\frac{\partial^{\beta}u(x_{i},y,t)}{\partial
|x|^{\beta}}=-\Delta x^{-\beta}\sum_{r=-1}^{1}\varrho_{r}^{(\beta)}\sum_{k=-\infty}^{\infty}\omega_{k}^{(\beta)}u(x_{i-k-r},y,t)
\end{eqnarray*}
\begin{eqnarray} \label{math201804010426}
+\mathcal{O}(\Delta x^{4}),\quad 0< \beta <1 ,
\end{eqnarray}

\begin{eqnarray*} \label{math201804010427}
\frac{\partial^{\mu}u(x,y_{j},t)}{\partial
|y|^{\mu}}=-\Delta y^{-\mu}\sum_{r=-1}^{1}\varrho_{r}^{(\mu)}\sum_{k=-\infty}^{\infty}\omega_{k}^{(\mu)}u(x,y_{j-k-r},t)
\end{eqnarray*}
\begin{eqnarray} \label{math201804010427}
+\mathcal{O}(\Delta y^{4}),\quad 1< \mu \leq2 ,
\end{eqnarray}
and
\begin{eqnarray*}\label{math201804010428}
\frac{\partial^{\nu}u(x,y_{j},t)}{\partial
|y|^{\nu}}=-\Delta y^{-\nu}\sum_{r=-1}^{1}\varrho_{r}^{(\nu)}\sum_{k=-\infty}^{\infty}\omega_{k}^{(\nu)}u(x,y_{j-k-r},t)
\end{eqnarray*}
\begin{eqnarray}\label{math201804010428}
+\mathcal{O}(\Delta y^{4}),\quad 0< \nu <1.
\end{eqnarray}

Let $u_{i,j}(t)=u(x_{i},y_{j},t)$, for $i=1,2,...,m_{1}-1$ and $j=1,2,...,m_{2}-1$, by considering zero boundary conditions, the two dimensional RSFADE (\ref{math201803252034}) can be cast into the following system of time ordinary differential equations by considering formulas (\ref{math201804010425})-(\ref{math201804010428}) based on mesh sizes in the spatial direction.
\begin{eqnarray*}
\frac{\partial u_{i,j}(t)}{\partial t}=-d_{\alpha}\Delta x^{-\alpha}\sum_{r=-1}^{1}\varrho_{r}^{(\alpha)}\sum_{k=-m_{1}+i}^{i}\omega_{k}^{(\alpha)}u_{i-k-r,j}(t)-c_{\beta}\Delta x^{-\beta}\sum_{r=-1}^{1}\varrho_{r}^{(\beta)}\times
\end{eqnarray*}
\begin{eqnarray*}
\sum_{k=-m_{1}+i}^{i}\omega_{k}^{(\beta)}u_{i-k-r,j}(t)-d_{\mu}\Delta y^{-\mu}\sum_{r=-1}^{1}\varrho_{r}^{(\mu)}\sum_{k=-m_{2}+j}^{j}\omega_{k}^{(\mu)}u_{i,j-k-r}(t)
\end{eqnarray*}
\begin{eqnarray}
-c_{\nu}\Delta y^{-\nu}\sum_{r=-1}^{1}\varrho_{r}^{(\nu)}\sum_{k=-m_{2}+j}^{j}\omega_{k}^{(\nu)}u_{i,j-k-r}(t)+s_{i,j}(t),
\end{eqnarray}
and by substituting $i-k-r=l$ and $j-k-r=z$, we have
\begin{eqnarray*}
\frac{\partial u_{i,j}(t)}{\partial t}=-d_{\alpha}(\sum_{r=-1}^{1}\varrho_{r}^{(\alpha)}\sum_{l=-r+m_{1}}^{-r}\frac{\omega_{i-l-r}^{(\alpha)}}{\Delta x^{\alpha}})u_{l,j}(t)-c_{\beta}(\sum_{r=-1}^{1}\varrho_{r}^{(\beta)}
\end{eqnarray*}
\begin{eqnarray*}
\times\sum_{l=-r+m_{1}}^{-r}\frac{\omega_{i-l-r}^{(\beta)}}{\Delta x^{\beta}})u_{l,j}(t)-d_{\mu}(\sum_{r=-1}^{1}\varrho_{r}^{(\mu)}\sum_{l=-r+m_{2}}^{-r}\frac{\omega_{i-l-r}^{(\mu)}}{\Delta x^{\mu}})u_{i,z}(t)
\end{eqnarray*}
\begin{eqnarray}
-c_{\nu}(\sum_{r=-1}^{1}\varrho_{r}^{(\nu)}\sum_{z=-r+m_{2}}^{-r}\frac{\omega_{j-z-r}^{(\nu)}}{\Delta y^{\nu}})u_{i,z}(t)+s_{i,j}(t),
\end{eqnarray}
i.e.,
\begin{eqnarray*}
\frac{\partial u_{i,j}(t)}{\partial t}=-(d_{\alpha}\sum_{r=-1}^{1}\varrho_{r}^{(\alpha)}\sum_{l=-r+m_{1}}^{-r}\frac{\omega_{i-l-r}^{(\alpha)}}{\Delta x^{\alpha}}+c_{\beta}\sum_{r=-1}^{1}\varrho_{r}^{(\beta)}
\end{eqnarray*}
\begin{eqnarray*}
\times\sum_{l=-r+m_{1}}^{-r}\frac{\omega_{i-l-r}^{(\beta)}}{\Delta x^{\beta}})u_{l,j}(t)-(d_{\mu}\sum_{r=-1}^{1}\varrho_{r}^{(\mu)}\sum_{l=-r+m_{2}}^{-r}\frac{\omega_{i-l-r}^{(\mu)}}{\Delta x^{\mu}}
\end{eqnarray*}
\begin{eqnarray}\label{Current3:4}
+c_{\nu}\sum_{r=-1}^{1}\varrho_{r}^{(\nu)}\sum_{z=-r+m_{2}}^{-r}\frac{\omega_{j-z-r}^{(\nu)}}{\Delta y^{\nu}})u_{i,z}(t)+s_{i,j}(t).
\end{eqnarray}
In order to display simple of the four operators in (\ref{Current3:4}) we show just briefly,
\begin{eqnarray*}
\mathfrak{D}_{\alpha,x}u_{i,j}^{n}=\sum_{r=-1}^{1}\varrho_{r}^{(\alpha)}\sum_{l=-r+m_{1}}^{-r}\frac{\omega_{i-l-r}^{(\alpha)}}{\Delta x^{\alpha}}u_{l,j}^{n},\quad \mathfrak{D}_{\beta,x}u_{i,j}^{n}=\sum_{r=-1}^{1}\varrho_{r}^{(\beta)}\sum_{l=-r+m_{1}}^{-r}\frac{\omega_{i-l-r}^{(\beta)}}{\Delta x^{\beta}}u_{l,j}^{n},
\end{eqnarray*}
\begin{eqnarray*}
\mathfrak{D}_{\mu,y}u_{i,j}^{n}=\sum_{r=-1}^{1}\varrho_{r}^{(\mu)}\sum_{l=-r+m_{2}}^{-r}\frac{\omega_{i-l-r}^{(\mu)}}{\Delta x^{\mu}}u_{i,z}^{n}\quad \mbox{and}\quad \mathfrak{D}_{\nu,y}u_{i,j}^{n}=\sum_{r=-1}^{1}\varrho_{r}^{(\nu)}\sum_{z=-r+m_{2}}^{-r}\frac{\omega_{j-z-r}^{(\nu)}}{\Delta y^{\nu}}u_{i,z}^{n}.
\end{eqnarray*}
The two-dimensional Meerschaert--Tadjeran finite difference method can be formulated as follows:
\begin{equation*}
\frac{u_{i,j}^{n+1}-u_{i,j}^{n}}{\Delta t}=-(d_{\alpha}\mathfrak{D}_{\alpha,x}+c_{\beta}\mathfrak{D}_{\beta,x}+d_{\mu}\mathfrak{D}_{\mu,y}+c_{\nu}\mathfrak{D}_{\nu,y})u_{i,j}^{n}+s_{i,j}^{n}
\end{equation*}
A Crank-Nicolson type finite difference equation for the two-dimensional Riesz space fractional advection-dispersion equation (\ref{math201803252034}) may be obtained by substituting the modified fractional centered difference estimates into the differential equation centered at time $t_{n+\frac{1}{2}}=\frac{1}{2}(t_{n+1}+t_{n})$ may be expressed as follows \cite{Tadjeran2007813}:
\begin{equation*}\label{math201804171858}
\frac{u_{i,j}^{n+1}-u_{i,j}^{n}}{\Delta t}=-(d_{\alpha}\mathfrak{D}_{\alpha,x}+c_{\beta}\mathfrak{D}_{\beta,x}+d_{\mu}\mathfrak{D}_{\mu,y}+c_{\nu}\mathfrak{D}_{\nu,y})\frac{u_{i,j}^{n+1}+u_{i,j}^{n}}{2}
\end{equation*}
\begin{equation}\label{math201804171858}
+\frac{s_{i,j}^{n+1}+s_{i,j}^{n}}{2}
\end{equation}
Now Eq. (\ref{math201804171858}) can be re-arranged and rewritten as this
\begin{equation*}
[\mathcal{I}+\frac{\Delta t}{2}(d_{\alpha}\mathfrak{D}_{\alpha,x}+c_{\beta}\mathfrak{D}_{\beta,x}+d_{\mu}\mathfrak{D}_{\mu,y}+c_{\nu}\mathfrak{D}_{\nu,y})]u_{i,j}^{n+1}=
\end{equation*}
\begin{equation}\label{math201804171932}
[\mathcal{I}-\frac{\Delta t}{2}(d_{\alpha}\mathfrak{D}_{\alpha,x}+c_{\beta}\mathfrak{D}_{\beta,x}+d_{\mu}\mathfrak{D}_{\mu,y}+c_{\nu}\mathfrak{D}_{\nu,y})]u_{i,j}^{n}+s_{i,j}^{n+\frac{1}{2}}\Delta t,
\end{equation}
where $\mathcal{I}$ is an identify operator (depends on the type of spatial operator) and $s_{i,j}^{n+\frac{1}{2}}=\frac{s_{i,j}^{n+1}+s_{i,j}^{n}}{2}$.\\
We note that the finite difference method (\ref{math201804171932}) can be rewritten as the following directional splitting factorization form by ignoring the high order perturbation terms
\begin{equation*}
[\mathcal{I}+\frac{\Delta t}{2}(d_{\alpha}\mathfrak{D}_{\alpha,x}+c_{\beta}\mathfrak{D}_{\beta,x})][\mathcal{I}+\frac{\Delta t}{2}
(d_{\mu}\mathfrak{D}_{\mu,y}+c_{\nu}\mathfrak{D}_{\nu,y})]u_{i,j}^{n+1}=
\end{equation*}
\begin{equation}\label{math201803181534}
[\mathcal{I}-\frac{\Delta t}{2}(d_{\alpha}\mathfrak{D}_{\alpha,x}+c_{\beta}\mathfrak{D}_{\beta,x})][\mathcal{I}-\frac{\Delta t}{2}
(d_{\mu}\mathfrak{D}_{\mu,y}+c_{\nu}\mathfrak{D}_{\nu,y})]u_{i,j}^{n}+s_{i,j}^{n+\frac{1}{2}}\Delta t
\end{equation}
To implement the practical reiterate scheme, by following the Peaceman-Rachford strategy, we can split the previous equation in two, by introducing an intermediate variable $u^ {*} $, which represents a solution computed at an intermediate time. Therefore, we obtain a type of Peaceman-Rachford ADI \cite{Peaceman195541},
\begin{equation}\label{math201803181535}
[\mathcal{I}+\frac{\Delta t}{2}(d_{\alpha}\mathfrak{D}_{\alpha,x}+c_{\beta}\mathfrak{D}_{\beta,x})]u_{i,j}^{*}=[\mathcal{I}-\frac{\Delta t}{2}
(d_{\mu}\mathfrak{D}_{\mu,y}+c_{\nu}\mathfrak{D}_{\nu,y})]u_{i,j}^{n}+s_{i,j}^{n+\frac{1}{2}}\Delta t
\end{equation}
\begin{equation*}\label{math201803181536}
[\mathcal{I}+\frac{\Delta t}{2}
(d_{\mu}\mathfrak{D}_{\mu,y}+c_{\nu}\mathfrak{D}_{\nu,y})]u_{i,j}^{n+1}=[\mathcal{I}-\frac{\Delta t}{2}(d_{\alpha}\mathfrak{D}_{\alpha,x}+c_{\beta}\mathfrak{D}_{\beta,x})]u_{i,j}^{*}
\end{equation*}
\begin{equation}\label{math201803181536}
+s_{i,j}^{n+\frac{1}{2}}\Delta t
\end{equation}
In order to solve (\ref{math201803252034})--(\ref{math201803252036}) the following algorithm can be used.\\
\noindent\textbf{Algorithm}\\
\textbf{Step (1)} First solve on each fixed horizontal slice $y=y_{k}$ $(k=1,2,...,m_{2}-1)$, a set of $m_{1}-1$ equations at the points $x_{i}$, $i=1,2,...,m_{1}-1$ defined by (\ref{math201803181535}) to obtain the middle solution slice $u_{i,k}^{*}$.\\
\textbf{Step (2)} Next alternating the spatial direction, and for each $x = x_{k}$ $(k=1,2,...,m_{1}-1)$
solving a set of $m_{2}-1$ equations defined by (\ref{math201803181536}) at the points $y_{j}$, $j = 1,2,...,m_{2}-1$,
to get $u_{i,j}^{n+1}$.

\section{Analysis of the ADI Crank-Nicolson scheme}

In this section, we bring essential theorems for convergence of the proposed method.

\begin{thm}\label{math201804192159}
The ADI-Crank Nicolson discretization for (\ref{math201803252034})-(\ref{math201803252036}) defined by (\ref{math201803181534}) is consistent, with a truncation error of the order $\mathcal{O}(\Delta x^{4})+\mathcal{O}(\Delta y^{4})+\mathcal{O}(\Delta t^2)$.
\end{thm}
\begin{proof}
First reminder that, as in the Crank-Nicolson method, for discretization the first derivative in the temporal direction $\partial u(x_{i},y_{j},t_{n+1/2})/\partial t$ gives the second-order accuracy via the centered divided difference $[u(x_{i},y_{j},t_{n+1})-u(x_{i},y_{j},t_{n})]/\Delta t$.

But in the spatial directions, we may employ Theorem \ref{math201802221502} to express accuracy
\begin{eqnarray*}\label{math201802221503}
\frac{\partial^{\alpha}u(x_{i},y_{j},t_{n})}{\partial
|x|^{\alpha}}=-\mathfrak{D}_{\alpha,x}u_{i,j}^{n}+\mathcal{O}(\Delta x^{4}),\quad \frac{\partial^{\beta}u(x_{i},y_{j},t_{n})}{\partial
|x|^{\beta}}=-\mathfrak{D}_{\beta,x}u_{i,j}^{n}+\mathcal{O}(\Delta x^{4}),
\end{eqnarray*}
\begin{eqnarray*}\label{math201802221503}
\frac{\partial^{\mu}u(x_{i},y_{j},t_{n})}{\partial
|y|^{\mu}}=-\mathfrak{D}_{\mu,y}u_{i,j}^{n}+\mathcal{O}(\Delta y^{4})\quad \mbox{and} \quad \frac{\partial^{\nu}u(x_{i},y_{j},t_{n})}{\partial
|y|^{\nu}}=-\mathfrak{D}_{\nu,y}u_{i,j}^{n}+\mathcal{O}(\Delta y^{4}).
\end{eqnarray*}

According to the above description, the local truncation error of the modified fractional centered difference Crank-Nicolson method presented in (\ref{math201804171858}) to approximate Eq. (\ref{math201803252034}), which is centered in both space and time, are forth order accurate in space and second order accurate in time.\\
By multiplying the sides of the relation (\ref{math201804171858}) in $\Delta t$ and by sorting terms based on algorithmic form in the temporal direction, the local truncation error of (\ref{math201804171932}) will be in order $\mathcal{O}(\Delta x^{4}\Delta t)+\mathcal{O}(\Delta y^{4}\Delta t)+\mathcal{O}(\Delta t^{3})$.\\
The interactions that apply to formula (\ref{math201803181534}) from formula (\ref{math201804171932}) will not affect the accuracy of the method. Since it is supposed to repeat the algorithm (\ref{math201803181534}) for the number of mesh points. Therefore, the ADI-CN finite difference equations defined by (\ref{math201803181534}) have a truncation error of the order $\mathcal{O}(\Delta x^{4})+\mathcal{O}(\Delta y^{4})+\mathcal{O}(\Delta t^2)$.
\end{proof}

\begin{thm}\label{math201804192200}
The iterative scheme defined by (\ref{math201803181534}) to solve the two dimensional RSFADE (\ref{math201803252034})-(\ref{math201803252036}) is unconditionally stable.
\end{thm}
\begin{proof}
Let $M_{x}$ and $M_{y}$ be matrix forms of the $\frac{\Delta t}{2}(d_{\alpha}\mathfrak{D}_{\alpha,x}+c_{\beta}\mathfrak{D}_{\beta,x})$ and $\frac{\Delta t}{2}
(d_{\mu}\mathfrak{D}_{\mu,y}+c_{\nu}\mathfrak{D}_{\nu,y})$, respectively. i.e,
\begin{equation*}
M_{x}=\frac{\Delta t}{2}m(d_{\alpha}\mathfrak{D}_{\alpha,x}+c_{\beta}\mathfrak{D}_{\beta,x})\quad \mbox{and} \quad M_{y}=\frac{\Delta t}{2}
m(d_{\mu}\mathfrak{D}_{\mu,y}+c_{\nu}\mathfrak{D}_{\nu,y})
\end{equation*}

According to Lemma \ref{math201803072254}, the matrices of $\mathfrak{D}_{\alpha,x}$, $\mathfrak{D}_{\beta,x}$, $\mathfrak{D}_{\mu,y}$ and $\mathfrak{D}_{\nu,y}$ operators  are all symmetric positive definite. Therefore the matrices $M_{x}$ and $M_{y}$ are symmetric positive definite. And the matrices used in the recurrence relation corresponded in formula (28) are $(I+M_{x})^{-1}(I-M_{x})$ and $(I+M_{x})^{-1}$ in the direction of the length axis and $(I+M_{y})^{-1}(I-M_{y})$ and $(I+M_{y})^{-1}$ in the direction of the latitude axis. Since the two matrices $M_{x}$ and $M_{y}$ are symmetric positive definite, the eigenvalues of these two matrices are positive. So the eigenvalues of the two matrices $(I+M_{x})^{-1}$ and $(I+M_{y})^{-1}$ are positive and the eigenvalues of the matrices $(I+M_{x})^{-1}(I-M_{x})$ and $(I+M_{y})^{-1}(I-M_{y})$ be inner the unit circle. In addition, it should be noted that two matrices $M_ {x} $ and $M_ {y} $ commute (are commutative) to the matrix multiplication. Therefore the iterative scheme defined by (\ref{math201803181534}) is unconditionally stable.
\end{proof}

Let us now recall the following Lax results that allows the authors to obtain fruitful considerations regarding the convergence.
\begin{thm}\label{math201804200911} \textbf{Lax Equivalence Theorem} \cite{Thomas2013}
A consistent, two level difference scheme for a well-posed linear initial-value problem is convergent
if and only if it is stable.
\end{thm}
\begin{thm}\label{math201804200914} \textbf{Lax Theorem} \cite{Thomas2013}
If a two-level difference scheme
\begin{equation*}
\textbf{u}^{n+1}=Q\textbf{u}^{n}+\Delta t \textbf{G}^{n}
\end{equation*}
is accurate of order $(p, q)$ in the norm $\parallel.\parallel$ to a well-posed linear initial-value
problem and is stable with respect to the norm $\parallel.\parallel$, then it is convergent of
order $(p, q)$ with respect to the norm $\parallel.\parallel$.
\end{thm}

The proposed method by referring to Theorem \ref{math201804192159} is either consistently or accurate in related order  and based on Theorem \ref{math201804192200} is unconditionally stable. Convergence of the method yields from Lax equivalence theorem and it follows from Lax Theorem that the convergence of numerical method is of order $(4,2)$ in spatial and temporal directions, respectively.

\section{Numerical results}

In this section, we report on numerical experiments for (\ref{math201803252034})-(\ref{math201803252036}) with the known exact solution to prove the correctness of our theoretical analysis expressed in the previous sections.

We use the same spacing $h$ in each direction, $\Delta x=\Delta y=h$. We present the error in $L^{2}$ norm
\begin{equation*}
e(h,k)=\parallel U^{N}-u^{N}\parallel,
\end{equation*}
 and the convergence rates in the temporal and spatial directions determined by the following
formulas
\begin{equation*}
\gamma_{1}\approx \frac{log(e_{m}/e_{m+1})}{log(h_{m}/h_{m+1})}, \quad \gamma_{2}\approx \frac{log(e_{n}/e_{n+1})}{log(k_{n}/k_{n+1})},
\end{equation*}
respectively, where the step size $h_{m}=Length/m$ , $k_{n}=T_{end}/n$ and $e_{m}$ , $e_{n}$ is the norm of the error with $h=h_{m}$ , $k=k_{n}$ respectively.

\begin{example}
\label{exp:1}
We consider the following two dimensional RSFADE with the initial and homogeneous Dirichlet boundary conditions:
\begin{eqnarray*}
\frac{\partial u(x,y,t)}{\partial
t}=d_{\alpha}\frac{\partial^{\alpha}u(x,y,t)}{\partial
|x|^{\alpha}}+c_{\beta}\frac{\partial^{\beta}u(x,y,t)}{\partial
|x|^{\beta}}+d_{\mu}\frac{\partial^{\mu}u(x,y,t)}{\partial
|y|^{\mu}}+c_{\nu}\frac{\partial^{\nu}u(x,y,t)}{\partial
|y|^{\nu}}
\end{eqnarray*}
\begin{eqnarray*}
+s(x,y,t),\quad 0<t<\pi ,\quad 0<x,y<1,
\end{eqnarray*}
\begin{eqnarray*}
u(x,y,0)=0, \quad 0\leq x,y \leq 1,
\end{eqnarray*}
\begin{eqnarray*}
u(0,y,t)=u(1,y,t)=u(x,0,t)=u(x,1,t)=0, \quad 0\leq t \leq \pi,\quad 0\leq x,y \leq 1,
\end{eqnarray*}
with source term
\begin{equation*}
s(x,y,t)=d_{\alpha}y^{2}(1-y)^{2}\frac{\sin(\pi t)}{\cos(\frac{\pi\alpha}{2})}\Phi(\alpha,x)
+c_{\beta}y^{2}(1-y)^{2}\frac{\sin(\pi t)}{\cos(\frac{\pi\beta}{2})}\Phi(\beta,x)
\end{equation*}
\begin{equation*}
+d_{\mu}x^{2}(1-x)^{2}\frac{\sin(\pi t)}{\cos(\frac{\pi\mu}{2})}\Phi(\mu,y)+c_{\nu}x^{2}(1-x)^{2}\frac{\sin(\pi t)}{\cos(\frac{\pi\nu}{2})}\Phi(\nu,y)
\end{equation*}
\begin{equation*}
+\pi x^{2}y^{2}(1-x)^{2}(1-y)^{2}\cos(\pi t)
\end{equation*}

where $\alpha=1.8$, $\beta=0.9$, $\mu=1.6$, $\nu=0.7$, $d_{\alpha}=d_{\mu}=0.25$, $c_{\beta}=c_{\nu}=0.05$ and $\Phi(\gamma,z)=\frac{z^{2-\gamma}+(1-z)^{2-\gamma}}{\Gamma(3-\gamma)}-\frac{6(z^{3-\gamma}+
(1-z)^{3-\gamma})}{\Gamma(4-\gamma)}+\frac{12(z^{4-\gamma}+(1-z)^{4-\gamma})}{\Gamma(5-\gamma)}$. The exact solution of the problem is given by $u(x,y,t)=x^{2}y^{2}(1-x)^{2}(1-y)^{2}\sin(\pi t)$.

The maximum absolute errors and their estimated convergence rates to solve the two dimensional RSFADE with the initial value and homogeneous Dirichlet boundary conditions approximated by the modified Crank-Nicolson ADI method are shown in Tables \ref{tab:1} and \ref{tab:2}. We observe that as the step sizes reduce, the absolute errors decrease as well, which shows that the numerical solutions is coincident with the exact solution. From the numerics in Tables \ref{tab:1} and \ref{tab:2}, we can see that the convergence rate for space is fourth order and the convergence rate for time is second order. The order of convergence is evaluated numerically, which demonstrates the theoretical results.
\begin{table}[ht!]
\caption{The maximum errors and convergence rates for the modified Crank-Nicolson ADI method for solving two dimensional RSFADE with halved spatial step sizes and $\mathit{k}_{t}= 0.001$} 
\label{tab:1} 
\centering 
\begin{tabular}{l l l l l } 
\hline
             & Maximum        & Estimated  \\
$\mathit{h}_{x}=\mathit{h}_{y}$& Absolute Error & Convergence Rate  \\
\hline \hline
$0.10000$ & $3.19826e-003$  &    -      \\
$0.05000$ & $2.61740e-004$  & $3.61108$ \\
$0.02500$ & $1.90572e-005$  & $3.77973$ \\
$0.01250$ & $1.33477e-006$  & $3.83567$ \\
$0.00625$ & $8.92357e-008$  & $3.90283$ \\
\hline \hline
\end{tabular}
\end{table}

\begin{table}[ht!]
\caption{The maximum errors and convergence rates for the modified Crank-Nicolson ADI method for solving two dimensional RSFADE with halved temporal step sizes and $\mathit{h}_{x}=\mathit{h}_{y}= 0.001$} 
\label{tab:2} 
\centering 
\begin{tabular}{l l l l l } 
\hline
             & Maximum        & Estimated  \\
$\mathit{k}_{t}$& Absolute Error & Convergence Rate  \\
\hline \hline
$0.10000$ & $3.77425e-003$  &    -      \\
$0.05000$ & $1.20417e-003$  & $1.64815$ \\
$0.02500$ & $3.55418e-004$  & $1.76045$ \\
$0.01250$ & $1.02360e-004$  & $1.79586$ \\
$0.00625$ & $2.69907e-005$  & $1.92312$ \\
\hline \hline
\end{tabular}
\end{table}

\end{example}

\begin{example}
\label{exp:2}
We consider the following two dimensional RSFADE with the initial and homogeneous Dirichlet boundary conditions:
\begin{eqnarray*}
\frac{\partial u(x,y,t)}{\partial
t}=d_{\alpha}\frac{\partial^{\alpha}u(x,y,t)}{\partial
|x|^{\alpha}}+c_{\beta}\frac{\partial^{\beta}u(x,y,t)}{\partial
|x|^{\beta}}+d_{\mu}\frac{\partial^{\mu}u(x,y,t)}{\partial
|y|^{\mu}}+c_{\nu}\frac{\partial^{\nu}u(x,y,t)}{\partial
|y|^{\nu}}
\end{eqnarray*}
\begin{eqnarray*}
+s(x,y,t),\quad 0<t<2,\quad 0<x,y<\pi,
\end{eqnarray*}
\begin{eqnarray*}
u(x,y,0)=xy(\pi-x)(\pi-y), \quad 0\leq x,y \leq \pi,
\end{eqnarray*}
\begin{eqnarray*}
u(0,y,t)=u(\pi,y,t)=u(x,0,t)=u(x,\pi,t)=0, \quad 0\leq t \leq 2,\quad 0\leq x,y \leq \pi,
\end{eqnarray*}
with source function
\begin{equation*}
s(x,y,t)=\frac{d_{\alpha}y(\pi-y)e^{-t}}{2\cos(\frac{\pi\alpha}{2})}\Psi(\alpha,x)
+\frac{c_{\beta}y(\pi-y)e^{-t}}{2\cos(\frac{\pi\beta}{2})}\Psi(\beta,x)
\end{equation*}
\begin{equation*}\label{5.2}
+\frac{d_{\mu}x(\pi-x)e^{-t}}{2\cos(\frac{\pi\mu}{2})}\Psi(\mu,y)+\frac{c_{\nu}x(\pi-x)e^{-t}}
{2\cos(\frac{\pi\nu}{2})}\Psi(\nu,y)-xy(\pi-x)(\pi-y)e^{-t}
\end{equation*}
where $\alpha=1.8$, $\beta=0.7$, $\mu=1.6$, $\nu=0.5$, $d_{\alpha}=d_{\mu}=0.25$, $c_{\beta}=c_{\nu}=0.05$ and $\Psi(\gamma,z)=\frac{\pi(z^{1-\gamma}+(\pi-z)^{1-\gamma})}{\Gamma(2-\gamma)}-\frac{2(z^{2-\gamma}+(\pi-z)^{2-\gamma})}
{\Gamma(3-\gamma)}$. The corresponding exact solution is $u(x,y,t)=xy(\pi-x)(\pi-y)e^{-t}$.

The table \ref{tab:3} shows maximum absolute errors and related estimated convergence rates with different values for  $\mathit{h}_{x}=\mathit{h}_{y}$ as $0.1\pi$, $0.05\pi$, $0.025\pi$, $0.0125\pi$ and $0.00625\pi$ and fixed value $\mathit{k}_{t}=0.001$ whereas Table \ref{tab:4} presents them with different values for $\mathit{k}_{t}$ as $0.1$, $0.05$, $0.025$, $0.0125$ and $0.00625$ and fixed value $\mathit{h}_{x}=\mathit{h}_{y}=0.001\pi$. From Tables \ref{tab:3} and \ref{tab:4}, we find the experimental convergence orders are approximately four and two in spatial and temporal directions, respectively. The numerical Example \ref{exp:2} results are provided to show that the proposed approximation method is computationally efficient.
\begin{table}[ht!]
\caption{The maximum errors and convergence rates for the modified Crank-Nicolson ADI method for solving two dimensional RSFADE with halved spatial step sizes and $\mathit{k}_{t}= 0.001$} 
\label{tab:3} 
\centering 
\begin{tabular}{l l l l l } 
\hline
             & Maximum        & Estimated  \\
$\mathit{h}_{x}=\mathit{h}_{y}$& Absolute Error & Convergence Rate  \\
\hline \hline
$0.10000\pi$ & $3.26587e-004$  &    -      \\
$0.05000\pi$ & $2.60038e-005$  & $3.65067$ \\
$0.02500\pi$ & $1.81670e-006$  & $3.83933$ \\
$0.01250\pi$ & $1.26448e-007$  & $3.84471$ \\
$0.00625\pi$ & $8.13028e-009$  & $3.95909$ \\
\hline \hline
\end{tabular}
\end{table}

\begin{table}[ht!]
\caption{The maximum errors and convergence rates for the modified Crank-Nicolson ADI method for solving two dimensional RSFADE with halved temporal step sizes and $\mathit{h}_{x}=\mathit{h}_{y}= 0.001\pi$} 
\label{tab:4} 
\centering 
\begin{tabular}{l l l l l } 
\hline
             & Maximum        & Estimated  \\
$\mathit{k}_{t}$& Absolute Error & Convergence Rate  \\
\hline \hline
$0.10000$ & $4.79240e-003$  &    -      \\
$0.05000$ & $1.46420e-003$  & $1.71063$ \\
$0.02500$ & $4.21902e-004$  & $1.79513$ \\
$0.01250$ & $1.14998e-004$  & $1.87530$ \\
$0.00625$ & $2.95945e-005$  & $1.95821$ \\
\hline \hline
\end{tabular}
\end{table}

\end{example}

\section{Conclusions}

In this paper, a high-order Crank-Nicolson alternating direction implicit method to solve two-dimensional Riesz space fractional advection-dispersion equations was introduced. Consistency and stability of the proposed method are discussed and by using the Lax equivalence theorem and Lax theorem, convergence and its order are proved. Numerical tests are carried out to indicate the theoretical results and demonstrate the efficiency of the proposed method.

\section*{Acknowledgments}

This work of the first and third authors is partially supported by Grant-in-Aid from the University of Mohaghegh Ardabili, Ardabil, Iran.\\
The second author is supported by the Ministry of Education and Science of the Russian Federation (5-100 program of the Russian Ministry of Education).

\end{document}